\documentclass[english]{article}
\usepackage[T1]{fontenc}
\usepackage[latin9]{inputenc}
\usepackage{babel}
\usepackage{amsfonts,color}

\newtheorem{dfn}{Definition}[section]
\newtheorem{theorem}[dfn]{Theorem}
\newtheorem{prop}[dfn]{Proposition}
\newtheorem{remark}[dfn]{Remark}

\newtheorem{lem}[dfn]{Lemma}

\usepackage{amssymb} 
\usepackage{amsmath}

 \global\long\def\sbr#1{\left[ #1\right] }
 \global\long\def\cbr#1{\left\{  #1\right\}  }
 
 \global\long\def\rbr#1{\left(#1\right)}

 \global\long\def\R{\mathbb{R}}
  \global\long\def\N{\mathbb{N}}
  
  \global\long\def\Z{\mathbb{Z}}

 \global\long\def\dd#1{\textnormal{d}#1}

 \global\long\def\TTV#1#2#3{\text{TV}^{#3}\!\rbr{#1,#2}}
 
 \global\long\def\UTV#1#2#3{\text{UTV}^{#3}\!\rbr{#1,#2}}
 \global\long\def\DTV#1#2#3{\text{DTV}^{#3}\!\rbr{#1,#2}}

 \global\long\def\ra{\rightarrow}
 
 \global\long\def\ns{\infty}
 \global\long\def\8{\infty}

 \global\long\def\ns{\infty}
\global\long\def\Ucross#1#2#3{\text{u}^{#1}\!\rbr{#2,#3}}
\global\long\def\Dcross#1#2#3{\text{d}^{#1}\!\rbr{#2,#3}}
\global\long\def\cross#1#2#3{\text{n}^{#1}\!\rbr{#2,#3}}

\global\long\def\Ucrossemph#1#2#3{\emph{u}^{#1}\!\rbr{#2,#3}}
 \global\long\def\Dcrossemph#1#2#3{\emph{d}^{#1}\!\rbr{#2,#3}}
 \global\long\def\crossemph#1#2#3{\emph{n}^{#1}\!\rbr{#2,#3}}

\title{Local times of deterministic paths with finite variation}

\author{Darlington Hove, Farai J. Mhlanga,  Rafa{\l } M. {\L }ochowski and Phumlani L. Zondi}

\begin{document}

\maketitle

\begin{abstract}
In this note, we define the numbers of level crossings by a c{\`a}dl{\`a}g (RCLL) real function $x\colon [0,+\ns) \ra \R$ and, in analogy to the work of \cite{BertoinYor:2014} we prove that for $x$ with locally finite total variation these numbers are densities of relevant occupation measures associated with $x$. Next, depending on the regularity of $x$ and $f:\R \ra \R$, we derive change of variable formulas, which may be seen as analogous of the It\^o or Tanaka-Meyer formulas. Some of these formulas are present in \cite{BertoinYor:2014} but we also present some generalizations.
\end{abstract}

\section{Introduction}

The total variation of a function $x\colon [0,+\ns) \ra \R$ on the time interval $[s,t]$, $0 \le s < t < +\ns$, is defined as:
\begin{equation} \label{tvc_def}
\TTV x{\left[s,t\right]}{}:=\sup_{\pi \in \Pi(s,t)}\sum_{[u,v] \in \pi} \left| x_{v}-x_{u} \right|,
\end{equation}
where the supremum is taken over all finite partitions $\pi$ of the interval $[s,t]$, that is, finite sets of non overlapping (with disjoint interiors) sub intervals $[u,v]$ of $[s,t]$ such that $\bigcup_{[u,v] \in \pi}[u,v] = [s,t]$. The family of all such partitions is denoted by $\Pi(s,t)$.

$x$ is said to have locally finite total variation if $\TTV x{\left[s,t\right]}{}$ is finite for any $s, t$ satisfying $0 \le s < t < +\ns$. In \cite{BertoinYor:2014}, the authors introduce local times $\ell^z(t)$ for c{\`a}dl{\`a}g functions with locally finite variation and state change of variable formulas which may be viewed as analogs of similar, fundamental formulas known in stochastic calculus (It\^o formula or It\^o-Tanaka-Meyer formula). The introduced local times also allow us to generalize the Banach Indicatrix Theorem, see \cite{BanachIndicatrix:1925, Vitaliindicatrix:1926}, to the case of c{\`a}dl{\`a}g (RCLL)  functions. Other generalization of the Banach Indicatrix may be found in \cite{Lozinskiindicatrix:1948, Lozinskiindicatrix:1958, LochowskiColloquium:2017}.

In this short paper, we define numbers of level (up-, down-) crossings by a c{\`a}dl{\`a}g real function $x\colon [0,+\ns) \ra \R$, $\Ucross{z}{x}{\sbr{0,t}}$, $\Dcross{z}{x}{\sbr{0, t}}$, and prove that they may be also viewed as densities of relevant occupation measures associated with $x$. Our approach seems to be more uniform and we obtain, among others, a simpler (and more general) change of variable formula, namely
\begin{align*}
f\rbr{x_t} - f\rbr{x_0} =  & \int_{\R}\cbr{\Ucross{z}{x}{\sbr{0,t}} -\Dcross{z}{x}{\sbr{0, t}} }  g(z)\dd z  \nonumber \\
& = \int_{\R} \ell^z(t) g(z)\dd z + \sum_{0<s\le t,\Delta x_{s}\neq 0}\cbr{f\rbr{x_{s}}-f\rbr{x_{s-}} },
\end{align*}
for any locally Lipschitz function $f:\R \ra \R$ with derivative $g$. 
In order to define the summation over the jumps of a c{\`a}dl{\`a}g function, we apply the concept of summation over general sets. Let $I$ be a set, let $b\colon I\to \R$ be a real valued function and let $\mathcal{I}$ be the family of all finite subsets of $I$. Since $\mathcal{I}$ is directed when endowed with the order of inclusion~$\subseteq$, the summation over $I$ can be defined by 
\begin{equation}\label{Sum AC}
  \sum_{i\in I} b_i:=\lim_{\Gamma \in \mathcal{I}}\sum_{i\in \Gamma} b_i
\end{equation} 
as limit of a net, i.e., $\lim_{\Gamma \in \mathcal{I}}\sum_{i\in \Gamma} b_i=:l \in (-\infty,\infty)$ exists if, for any $\varepsilon >0$ there is  $\Gamma\in \mathcal{I}$ such that for all $\tilde \Gamma\in \mathcal{I}$ satisfying $\tilde \Gamma \supseteq \Gamma$ one has $ \sum_{i\in \tilde\Gamma} b_i\in \rbr{l - \varepsilon, l +\varepsilon}$. 

This paper is organised as follows: in Section 2, we introduce the necessary notations, definitions, and the main result of our paper. Section 3 provides an alternative proof to the main result. Finally, in Section 4, we present change of variable formulas.

\section{Definitions, notation and the main result}

By $\Z$ we denote the set of all integers, by $\N$ we denote the set of all positive integers and by $\N_0$ we denote the set $\N \cup \cbr{0}$. By $\mathbb{Q}$ we denote the set of all rational numbers. By $C([0,+\ns);\R)$ we denote the space of all continuous functions.
By $D([0,+\ns);\R)$ we denote the space of all c{\`a}dl{\`a}g (RCLL) functions $x\colon [0,+\ns) \ra \R$, that is $x \in D([0,+\ns);\R)$ if it is right-continuous at each $t\in [0,+\ns)$ and possesses left-limits at each $t\in (0,+\ns)$.
Let $V = V([0,+\ns);\R)$ be the subset of $D = D([0,+\ns);\R)$ of functions with locally finite total variation.

Let $y\in \R, c > 0$. In what follows, we introduce \emph{numbers of (up-, down-) crossings}  of the interval $[y-c/2,y+c/2]$ by $x \in D([0,+\ns);\R)$ during the time interval $[s, t]$.
We define $\sigma_{0}^{c}=s$ and for $n \in \N_0$
\[
\tau_{n}^{c}=\inf\left\{ u \in \sbr{ \sigma_{n}^{c}, t}: f(u)\ge y+c/2\right\} ,
\]
\[
\sigma_{n+1}^{c}=\inf\left\{ u \in \sbr{ \tau_{n}^{c}, t} : f(u) < y-c/2\right\},
\]
where we apply conventions: $\inf \emptyset = +\ns$, $\sbr{+\ns, t} = \emptyset$.
\begin{dfn}\label{defd}
The \emph{number of downcrossings by $x$ the interval $[y-c/2, y+c/2]$} during the time interval $[s,t]$ is defined as
\[
\Dcrossemph{y,c}{x}{[s,t]} :=\max\left\{ n:\sigma_{n}^{c}\leq t\right\} .
\]
\end{dfn}
\begin{dfn}\label{defu}
The \emph{number of upcrossings by $x$ the interval $[y-c/2, y+c/2]$} during the time interval $[s,t]$ is defined as
\[
\Ucrossemph{y,c}{x}{[s,t]} := \Dcrossemph{-y,c}{-x}{[s,t]}.
\]
\end{dfn}
\begin{dfn}\label{den}
The \emph{number of crossings by $x$ the interval $[y-c/2, y+c/2]$} during the time interval $[s,t]$ is defined as
\[
\crossemph{y,c}{x}{[s,t]}:= \Ucrossemph{y,c}{x}{[s,t]} + \Dcrossemph{y,c}{x}{[s,t]}.
\]
\end{dfn}

\begin{remark} \label{durepresentation}
 If $y+c/2$ is not a local maximum of $x$, then the number $\Dcrossemph{y,c}{x}{[s,t]}$ may be also defined as
  \begin{equation*}
    \Dcrossemph{y,c}{x}{[s,t]}:= \sup_{n\in \N} \sup_{s \leq s_1 <t_1 <\ldots <s_n<t_n \leq t; \, s_1,\ldots, s_n, t_1, \ldots, t_n \in \mathbb{Q}} \sum_{i=1}^n \mathrm{D}^{y,c}\rbr{x_{s_i}, x_{t_i}},
  \end{equation*}
  where
  \begin{equation*}
    \mathrm{D}^{y,c}(a, b):= \begin{cases}
                         1 & \text{if } a\ge y+{c}/{2} \text{ and } b <  y-{c}/{2},\\
                         0 & \text{if otherwise}.
                       \end{cases}
  \end{equation*}
  The numbers $\Ucrossemph{y,c}{x}{[s,t]}$ of upcrossings may be defined analogously.
 \end{remark}
  From Remark \ref{durepresentation}, the measurability of $y \mapsto \Dcross{y,c}{x}{[s,t]}$, $y \mapsto \Ucross{y,c}{x}{[s,t]}$ and $y \mapsto \cross{y,c}{x}{[s,t]}$ immediately follows.

Using numbers of interval (down-, up-) crossings we define {\em level (down-, up-) crossings}.
  \begin{dfn}\label{defd1} Let $y\in\R$. The number of times that the function $x$ {\em downcrosses} the level $y$ during the time interval $[s, t]$ is defined as
	  \[
	 	\Dcrossemph{y}{x}{[s,t]} = \lim_{c\to 0+}\Dcrossemph{y,c}{x}{[s,t]}\in \N_0 \cup\cbr{+\ns}.
	 \]
	 Analogously, we define the number of {\em upcrosses} as
	\[
	 	\Ucrossemph{y}{x}{[s,t]} = \lim_{c\to 0+}\Ucrossemph{y,c}{x}{[s,t]}\in \N_0 \cup\cbr{+\ns}.
	 \]
	 Finally, we define the number of times that the function $x$ {\em crosses} the level $y$ during the time interval $[s, t]$
	 \[
	 	\crossemph{y}{x}{[s,t]} = \lim_{c\to 0+}\crossemph{y,c}{x}{[s,t]}\in \N_0 \cup\cbr{+\ns}.
	 \]
\end{dfn}

Next, let us also define {\em level (down-, up-) crossings} via jumps. For $s \in (0, +\ns)$, $\Delta x_s$ denotes the jump of  $x$ at the time $s$, $\Delta x_s = x_{s} - x_{s-}$.
 \begin{dfn}\label{defd1} Let $y\in\R$. The number of times that the function $x$ {\em downcrosses} the level $y$ during the time interval $[s, t]$ via jumps is defined as
	  \[
	 	\Delta \Dcrossemph{y}{x}{[s,t]} =\text{\emph{Card}} \cbr{u \in (s,t]: \Delta x_u <0 \text{ and } y \in \rbr{x_u, x_{u-}} }.
	 \]
	 Analogously we define the number of {\em upcrosses via jumps} as
	\[
	 	\Delta \Ucrossemph{y}{x}{[s,t]} = \text{\emph{Card}} \cbr{u \in (s,t]: \Delta x_u > 0 \text{ and } y \in \rbr{x_{u-}, x_u} }.
	 \]
	 Finally, we define the number of times that the function $x$ {\em crosses via jumps} the level $y$ during the time interval $[s, t]$
	 \[
	 	 \Delta \crossemph{y}{x}{[s,t]} = \Delta \Dcrossemph{y}{x}{[s,t]} + \Delta \Ucrossemph{y}{x}{[s,t]} .
	 \]
\end{dfn}
\begin{remark} \label{Deltadurepresentation}
 The numbers $\Delta\Dcrossemph{y}{x}{[s,t]}$ may be aslo defined as
  \begin{align*}
  &  \Delta\Dcrossemph{y,c}{x}{[s,t]} \\
    &  := \sup_{c\in \mathbb{Q}, c>0} \inf_{m \in \N} \sup_{n\in \N} \sup_{s \leq s_1 <t_1 <\ldots <s_n<t_n \leq t; \, s_i, t_i \in \mathbb{Q}, t_i -s_i \le \frac{1}{m}, i=1,\ldots,n} \sum_{i=1}^n \mathrm{D}^{y,c}\rbr{x_{s_i}, x_{t_i}},
  \end{align*}
  where  the numbers $\mathrm{D}^{y,c}\rbr{a,b}$ were defined in Remark \ref{durepresentation}.

      The numbers $\Delta \Ucrossemph{y,c}{x}{[s,t]}$ of upcrossings may be defined analogously.
 \end{remark}
  From Remark \ref{Deltadurepresentation}, the measurability of $y \mapsto \Delta \Dcross{y,c}{x}{[s,t]}$, $y \mapsto \Delta \Ucross{y,c}{x}{[s,t]}$ and $y \mapsto \Delta\cross{y,c}{x}{[s,t]}$ immediately follows.

One of the main goals of this note is to prove a relationship between the just defined numbers of level crossings and the occupation measures associated with $x$ and with the continuous part of $x$ which we will denote by $x^{cont}$ and which is defined by the relation
\[
x(t) = x^{cont}(t) + \sum_{0<s \le t} \Delta x_s, \quad t\in [0, +\ns).
\]

Let us now define two companion quantities of the total variation:  -- \emph{positive} and \emph{negative total variation} or \emph{upward} and \emph{downward total variation}. They are defined respectively as
\begin{equation*}
\UTV x{\left[s,t\right]}{}:=\sup_{\pi \in \Pi(s,t)}\sum_{[u,v] \in \pi} \rbr{x_{v}-x_{u}}_+
\end{equation*}
and
\begin{equation*}
\DTV x{\left[s,t\right]}{}:=\sup_{\pi \in \Pi(s,t)}\sum_{[u,v] \in \pi} \rbr{x_{v}-x_{u}}_-,
\end{equation*}
where for a real $y$, $(y)_+ := \max(0,y)$, $(y)_- := \max(0,-y)$.

Let $|\mu|$ be defined as the unique measure on $[0,+\ns)$
which assigns to the set $\cbr{0}$ the value $0$ and assigns to the interval $(a,b]\subset [0,+\ns)$, such that $0\le a< b<+\ns$,
the value
\[
|\mu|(a,b]:=\TTV x{[0,b]}{}-\TTV x{[0,a]}{} = \TTV x{[a,b]}{}.
\]
In what follows we will denote the Lebesgue integrals $\int_{A}h\dd{|\mu|}$
by
\[
\int_{A}h_{s}|\dd{x_{s}}|.
\]
Two other associated measures are $\mu_+$ and  $\mu_-$ which are defined as the unique measures on $[0,+\ns)$
which assign to the set $\cbr{0}$ the value $0$ and assigns to the interval $(a,b]\subset [0,+\ns)$, such that $0\le a < b<+\ns$,
the values
\[
\mu_+(a,b]:=\UTV x{[0,b]}{}-\UTV x{[0,a]}{} = \UTV x{[a,b]}{}
\]
and
\[
\mu_-(a,b]:=\DTV x{[0,b]}{}-\DTV x{[0,a]}{} = \DTV x{[a,b]}{}
\]
respectively.
In what follows we will denote the Lebesgue integrals $\int_{A}h\dd{\mu_+}$ and $\int_{A}h\dd{\mu_-}$
by
\[
\int_{A}h_{s}\rbr{\dd x_{s}}_+ \text{ and }  \int_{A}h_{s}\rbr{\dd x_{s}}_-
\]
respectively.

Similarly as the integrals $\int_{A}h_{s}|\dd{x_{s}}|$, $\int_{A}h_{s}\rbr{\dd x_{s}}_+$ and $\int_{A}h_{s}\rbr{\dd x_{s}}_-$, one defines the integrals $\int_{A}h_{s}|\dd{x_{s}^{cont}}|$, $\int_{A}h_{s}\rbr{\dd x_{s}}_+$ and $\int_{A}h_{s}\rbr{\dd x_{s}}_-$. \\

Now we are ready to state the following.

\begin{theorem} \label{main} For any function $x$ from the
family $V$, if $t>0$ and $g:\R\ra\R$ is Borel-measurable and locally  bounded, then:
\begin{align}
& \int_{\R} \Ucrossemph{z}{x}{[0,t]} g(z)\dd z \nonumber \\ 
& =\int_{(0,t]}g\rbr{x_{s-}}\rbr{\dd{x_{s}}}_{+}+\sum_{0<s\le t,\Delta x_{s}>0}\int_{x_{s-}}^{x_{s}}\cbr{g(z)-g\rbr{x_{s-}}}\dd z\nonumber \\
 & =\int_{(0,t]}g\rbr{x_{s-}}\rbr{\dd{x_{s}}^{cont}}_{+}+\sum_{0<s\le t,\Delta x_{s}>0}\int_{x_{s-}}^{x_{s}}g(z)\dd z\nonumber \\
 & =\int_{(0,t]}g\rbr{x_{s}}\rbr{\dd{x_{s}}^{cont}}_{+}+\sum_{0<s\le t,\Delta x_{s}>0}\int_{x_{s-}}^{x_{s}}g(z)\dd z.\label{eq:BanInd1}
\end{align}
and
\begin{align}
& \int_{\R}\Dcrossemph{z}{x}{[0,t]} g(z)\dd z \nonumber \\
 & =\int_{(0,t]}g\rbr{x_{s-}}\rbr{\dd{x_{s}}}_{-}+\sum_{0<s\le t,\Delta x_{s}<0}\int_{x_{s}}^{x_{s-}}\cbr{g(z)-g\rbr{x_{s-}}}\dd z\nonumber \\
 & =\int_{(0,t]}g\rbr{x_{s-}}\rbr{\dd{x_{s}}^{cont}}_{-}+\sum_{0<s\le t,\Delta x_{s}<0}\int_{x_{s}}^{x_{s-}}g(z)\dd z\nonumber \\
 & =\int_{(0,t]}g\rbr{x_{s}}\rbr{\dd{x_{s}}^{cont}}_{-}+\sum_{0<s\le t,\Delta x_{s}<0}\int_{x_{s}}^{x_{s-}}g(z)\dd z.\label{eq:BanInd1-1}
\end{align}
\end{theorem}
\begin{remark} Companion formulas read as:
\begin{align}
& \int_{\R}\cbr{\Ucrossemph{z}{x}{[0,t]} -\Dcrossemph{z}{x}{[0,t]} }  g(z)\dd z \nonumber \\
& =\int_{(0,t]}g\rbr{x_{s-}}{\dd{x_{s}}}+\sum_{0<s\le t,\Delta x_{s}\neq 0}\int_{x_{s-}}^{x_{s}}\cbr{g(z)-g\rbr{x_{s-}}}\dd z\nonumber \\
 & =\int_{(0,t]}g\rbr{x_{s-}}{\dd{x_{s}}^{cont}}+\sum_{0<s\le t,\Delta x_{s} \neq 0}\int_{x_{s-}}^{x_{s}}g(z)\dd z\nonumber \\
 & =\int_{(0,t]}g\rbr{x_{s}}{\dd{x_{s}}^{cont}}+\sum_{0<s\le t,\Delta x_{s} \neq 0}\int_{x_{s-}}^{x_{s}}g(z)\dd z\label{eq:BanInd1-11}
\end{align}
and
\begin{align}
& \int_{\R} \crossemph{z}{x}{[0,t]} g(z)\dd z \nonumber \\
& =\int_{(0,t]}g\rbr{x_{s-}}\left|\dd{x_{s}}\right|+\sum_{0<s\le t,\Delta x_{s}\neq0}\int_{x_{s-}\wedge x_{s}}^{x_{s}\vee x_{s}}\cbr{g(z)-g\rbr{x_{s-}}}\dd z\nonumber \\
 & =\int_{(0,t]}g\rbr{x_{s-}}\left|\dd{x_{s}^{cont}}\right|+\sum_{0<s\le t,\Delta x_{s}\neq0}\int_{x_{s-}\wedge x_{s}}^{x_{s}\vee x_{s}}g(z)\dd z\nonumber \\
 & =\int_{(0,t]}g\rbr{x_{s}}\left|\dd{x_{s}^{cont}}\right|+\sum_{0<s\le t,\Delta x_{s}\neq0}\int_{x_{s-}\wedge x_{s}}^{x_{s}\vee x_{s}}g(z)\dd z.\label{eq:BanInd1-2}
\end{align}
\end{remark}
Theorem \ref{main} is almost an immediate consequence of \cite[Theorem 1]{BertoinYor:2014}. In \cite{BertoinYor:2014}, for a real $z$ the authors introduce positive times at which the function $x$ increases continuously through the level $z$ and times at which the function $x$ decreases continuously through the level $z$. They denote the sets of the former times by ${\cal I}(z)$  and the latter ${\cal D}(z)$. They are defined in the following way.
\begin{itemize}
\item
$t \in {\cal I}(z)$ iff
(i) $x$ is continuous at time $t$ and $x(t) = z$ and (ii) the sign of $x(t) - x(s)$ is the same as the sign of $t-s$ for all $s$ in some neighborhood of $t$.
 \item
$t \in {\cal D}(z)$ iff
(i) $x$ is continuous at time $t$ and $x(t) = z$ and (ii) the sign of $x(t) - x(s)$ is the same as the sign of $s-t$ for all $s$ in some neighborhood of $t$.
\end{itemize}
 For $t>0$ and $z \in \R$ Bertoin and Yor introduce also the number $n^z(t)$ (a.k.a. Banach's Indicatrix) defined as
 \[
 n^z(t) := \text{{Card}} \cbr{s \in [0,t]: x(s) = z }.
 \]
After introduction ${\cal I}(z)$, ${\cal D}(z)$ and $n^z(t)$ we are ready to prove Theorem \ref{main}. \\
{\bf Proof:} For $z \in \R$ and $t>0$,  let us denote
\[
r^z(t) := n^z(t)  - \text{{Card}} \rbr{ {\cal I}(z) \cap (0, t]} -  \text{{Card}} \rbr{ {\cal D}(z) \cap (0, t]} \ge 0.
\]
It is not difficult to see that for $z \in \R$ and $c,t>0$ the following estimates hold
\[
\text{{Card}} \rbr{ {\cal I}(z) \cap (0, t]} + \Delta \Ucross{z}{x}{[0,t]} \le \Ucross{z}{x}{[s,t]}
\]
and
\[
\Ucross{z,c}{x}{[s,t]} \le \text{{Card}} \rbr{ {\cal I}(z) \cap (0, t]} + \Delta \Ucross{z}{x}{[0,t]}  +r^z(t).
\]
Proceeding to the limit $c \ra 0+$ we obtain
\begin{align*}
\text{{Card}} \rbr{ {\cal I}(z) \cap (0, t]} + \Delta \Ucross{z}{x}{[0,t]} & \le \Ucross{z}{x}{[s,t]} \\
& \le \text{{Card}} \rbr{ {\cal I}(z) \cap (0, t]} + \Delta \Ucross{z}{x}{[0,t]}  +r^z(t)
\end{align*}
which is equivalent to
\begin{equation} \label{fiiirst}
\text{{Card}} \rbr{ {\cal I}(z) \cap (0, t]} \le \Ucross{y}{x}{[s,t]} -  \Delta \Ucross{z}{x}{[0,t]}\le \text{{Card}} \rbr{ {\cal I}(z) \cap (0, t]}  +r^z(t).
\end{equation}
By {[}Theorem 1, Bertoin Yor BLMS 2014{]}, we have
\begin{align}
& \int_{(0,t]}g\rbr{x_{s}}\rbr{\dd{x_{s}}^{cont}}_{+} = \frac{1}{2} \int_{(0,t]}g\rbr{x_{s}}\rbr{\dd{x_{s}}^{cont} + \left|\dd{x_{s}^{cont}}\right|} \nonumber \\
& = \frac{1}{2} \int_{\R}g\rbr{z} \cbr{\ell^z(t) + \lambda^z(t) } \dd z = \int_{\R}g\rbr{z}  \text{{Card}} \rbr{ {\cal I}(z) \cap (0, t]} \dd z, \label{seeecond}
\end{align}
where
\[
\ell^z(t) = \text{{Card}} \rbr{ {\cal I}(z) \cap (0, t]} - \text{{Card}} \rbr{ {\cal D}(z) \cap (0, t]},
\]
\[
\lambda^z(t)  = \text{{Card}} \rbr{ {\cal I}(z) \cap (0, t]} + \text{{Card}} \rbr{ {\cal D}(z) \cap (0, t]}.
\]
Next, we observe that
\begin{align} 
& \int_{\R}g\rbr{z} \cbr{ \Ucross{z}{x}{[0,t]} -  \Delta \Ucross{z}{x}{[0,t]}}  \dd z \nonumber \\ 
& = \int_{\R}g\rbr{z} \Ucross{z}{x}{[0,t]}  \dd z - \sum_{0<s\le t,\Delta x_{s}>0}\int_{x_{s-}}^{x_{s}}g(z)\dd z. \label{threee}
\end{align}
(All terms on the right side of \eqref{threee} are well defined since $x$ has locally finite variation and $g$ is locally  bounded.)
The last ingredient we need is the equality
\begin{equation} \label{nuuul}
 \int_{\R} g\rbr{z} r^z(t) \dd z = 0,
\end{equation}
which also follows from \cite[Theorem 1]{BertoinYor:2014}, namely from the equality of the measures $n^z(t) \dd z$ and $\lambda^z(t) \dd z$.
From \eqref{nuuul} and \eqref{seeecond}, we obtain that
\begin{equation} \label{fooour}
\int_{\R}g\rbr{z}  \cbr{ \text{{Card}} \rbr{ {\cal I}(z) \cap (0, t]} + r^z(t)}\dd z =  \int_{(0,t]}g\rbr{x_{s}}\rbr{\dd{x_{s}}^{cont}}_{+}.
\end{equation}
This, together with   \eqref{fiiirst}, \eqref{seeecond} and  \eqref{threee} also gives
\begin{equation} \label{fiiive}
 \int_{\R}g\rbr{z} \Ucross{z}{x}{[0,t]}  \dd z - \sum_{0<s\le t,\Delta x_{s}>0}\int_{x_{s-}}^{x_{s}}g(z)\dd z =  \int_{(0,t]}g\rbr{x_{s}}\rbr{\dd{x_{s}}^{cont}}_{+}.
\end{equation}
From the continuity of $x^{cont}$, we further obtain
\begin{align}
 \int_{(0,t]}g\rbr{x_{s}}\rbr{\dd{x_{s}}^{cont}}_{+} & = \int_{(0,t]}g\rbr{x_{s-}}\rbr{\dd{x_{s}}^{cont}}_{+} \nonumber \\
& =  \int_{(0,t]}g\rbr{x_{s-}}\rbr{\dd{x_{s}}}_{+} -  \sum_{0<s\le t,\Delta x_{s}>0}\int_{x_{s-}}^{x_{s}}g\rbr{x_{s-}}\dd z.  \label{siiix}
\end{align}
Relations from the last two displays give \eqref{eq:BanInd1}.

The proof of \eqref{eq:BanInd1-1} is similar.
\hfill $\square$

\section{Alternative proof of Theorem \ref{main}}

In the proof of Theorem \ref{main}, we used the main result of  \cite{BertoinYor:2014}, that is \cite[Theorem 1]{BertoinYor:2014}. To prove their theorem, Bertoin and Yor use change of variable formulas, which are analogs of the It\^o and Tanaka-M\'eyer formulas known from stochastic calculus, namely, if $f: \R \ra \R$ is a continuously differentiable function ($f \in C^1(\R, \R)$) then for $t>0$:
\begin{equation} \label{Itooo}
 f\rbr{x_t} - f\rbr{x_0} = \int_{0}^t f'\rbr{x_s} \dd x_s^{cont} + \sum_{0<s\le t,\Delta x_{s}>0} \rbr{ f\rbr{x_{s}} - f\rbr{x_{s-}}},
\end{equation}
\begin{equation} \label{ItoMeyyy}
 f\rbr{x_t} - f\rbr{x_0} = \int_{\R} f'\rbr{z} \ell^z(t) \dd z + \sum_{0<s\le t,\Delta x_{s}>0} \rbr{ f\rbr{x_{s}} - f\rbr{x_{s-}}}
\end{equation}
and
\begin{equation} \label{TanakaMeyyy}
{\bf 1}_{[z, +\ns)}\rbr{x_t}= {\bf 1}_{[z, +\ns)}\rbr{x_0} + \ell^z(t) +  \sum_{0<s\le t,\Delta x_{s}>0} \rbr{ {\bf 1}_{[z, +\ns)}\rbr{x_s} - {\bf 1}_{[z, +\ns)}\rbr{x_{s-}}}.
\end{equation}
Formula \eqref{TanakaMeyyy} holds under the condition that  $z$ is a \emph{simple level}, which means that for any $t>0$
\[
\text{Card}\cbr{s \in (0,t]: x_{s-}< z < x_s \text{ or } x_s < z < x_{s-} \text{ or } x_s = z} < +\ns
\]
and there is no jump of $x$ that starts or ends at the level $z$.

In what follows, we present an alternative proof, from which, depending on the regularity of $g$ or $x$, we extend formula \eqref{ItoMeyyy} to the case of less regular (not necessarily continuous) $g$.

Let $V^{0}$ denote the subset of $D$ (and $V)$ consisting of functions which are piecewise monotonic, that is, for any $t >0$ there exists a \emph{finite} sequence of intervals $\{I_i\}_{i=1}^N$, $N \in \N $, which are mutually disjoint, $\bigcup_{i=1}^N I_i = [0,t]$ and the function $x$ is monotone on any of $I_i$.

First, we will  prove \eqref{eq:BanInd1} and  \eqref{eq:BanInd1-1} for in the case when $x \in V^0$.

\begin{lem} \label{Ban_Viatli_extended}
Let $x\in V^0$, $t>0$ and $g:\R\mapsto \R$ be Borel--measurable and locally integrable. Then the equalities
 \eqref{eq:BanInd1} and  \eqref{eq:BanInd1-1} hold.
\end{lem}
{\bf Proof:}
 Let $\{I_i\}_{i=1}^N$, $N \in \N $ be a sequence of intervals  which are mutually disjoint, $\bigcup_{i=1}^N I_i = [0,t]$ and the function $x$ is monotone on any of $I_i$. Since $x$ is c\`adl\`ag, we may and will assume that these intervals, except the last one containing $t$, are of the form $I_i = \left[t_i, t_{i+1}\right)$, where $t_i < t_{i+1}$. Moreover, we will assume that they are the largest intervals possible on which $x$ is monotone. Thus, if for some $i>1$, $x$ is non-decreasing on $ \left[t_{i-1}, t_{i}\right)$ and non-decreasing on $ \left[t_i, t_{i+1}\right)$ then there need to be a downward jump at the time $t_i$. Similarly, if for some $i>1$, $x$ is non-increasing on $ \left[t_{i-1}, t_{i}\right)$ and non-increasing on $ \left[t_i, t_{i+1}\right)$  then there need to be an upward jump at the time $t_i$. Let $S$ denote the set of such times $s$.

Let us define $R = \bigcup_{i=1}^N \cbr{x_{t_i}, x_{t_{i+1}-}, x_{t_{i+1}}}$. Given a level $z \in \R \setminus R$ with any interval $I_i$ or any $s \in S$, we can have at most one associated upcrossing. This happens when $x_{t_i} < x_{t_{i+1}}$ or $x_{s-} < x_s$ and
	\[
		z\in \rbr{x_{t_i}, x_{t_{i+1}}} \text{ or } z\in \rbr{x_{s-}, x_s}.
	\]
	Clearly, neither of intervals where function is decreasing induces any upcrossing. Let $J \subset \cbr{1,2,\ldots, N}$ be the subset of such indices $i\in \cbr{1,2,\ldots, N}$ that the function is increasing on $\left[ {t_i}, {t_{i+1}} \right)$ and $T$ be the subset of $s \in S$ such that $x$ has upward jump at $s$.  Then we have
	\[
		\Ucross{z}{x}{[0,t]} = \sum_{i\in J} \mathbf{1}_{ \rbr{x_{t_i}, x_{t_{i+1}}}  } (z) + \sum_{s\in T} \mathbf{1}_{ \rbr{x_{s-}, x_s} } (z) .
	\]
Consequently
	\[
		\int_{\R} g(z)\Ucross{z}{x}{t} \dd{z} = \int_{\R \setminus R} g(z)\Ucross{z}{x}{t} \dd{z} = \sum_{i\in J} \int_{x_{t_i}}^{x_{t_{i+1}}} g(z) \dd{z} + \sum_{s\in T} \int_{x_{s-}}^{x_{s}} g(z) \dd{z}.
	\]
	We are now to deal with these integrals. To this end, we apply the classical idea of opening temporal windows at times of jumps. We may impose that the sum of the lengths of these windows is finite and consider interpolated continuous $\tilde x$. We have
\[
\int_{x_{t_i}}^{x_{t_{i+1}}} g(z) \dd{z} = \int_{\tilde{x}_u}^{\tilde{x}_{u'}}g(z)\dd z \text{ and }
\int_{x_{s-}}^{x_{s}} g(z) \dd{z} = \int_{\tilde{x}_v}^{\tilde{x}_{v'}}g(z)\dd z
\]
for properly defined $u,u', v, v'$. Then we apply the change of variable for the Riemann-Stieltjes integral
\[
\int_{\tilde{x}_u}^{\tilde{x}_{u'}}g(z)\dd z  = \int_{u}^{u'} g{\rbr{\tilde{x}_s}} \dd \tilde{x}_{s} = \int_{u}^{u'} g{\rbr{\tilde{x}_s}} \rbr{\dd{\tilde{x}_{s}}}_{+}
 \]
and
\[
\int_{\tilde{x}_v}^{\tilde{x}_{v'}}g(z)\dd z  = \int_{v}^{v'} g{\rbr{\tilde{x}_s}} \dd \tilde{x}_{s}
 = \int_{v}^{v'} g{\rbr{\tilde{x}_s}}\rbr{\dd{\tilde{x}_{s}}}_{+}.
\]
Clearly, for properly defined $\tilde{t}$ we have
\begin{equation} \label{change_of_var}
\int_{\R} g(z)\Ucross{z}{x}{t} \dd{z} = \int_{0}^{\tilde{t}}  g{\rbr{\tilde{x}_s}} \rbr{\dd{\tilde{x}_{s}}}_{+}.
\end{equation}
Let us consider the decomposition of the measure $\text{UTV}(x,\dd s) = \rbr{\dd x_s}_+$ into the continuous part $\mu^{cont}$ and the atomic part. Moreover, let $W$ be the set of the temporal windows. We obtain
\[
	 \int_{0}^{\tilde{t}} \mathbf{1}_{s\notin W} g{\rbr{\tilde{x}_s}} \rbr{\dd{\tilde{x}_{s}}}_{+} = \int_0^t g(x_{s-}) \dd \mu^{cont}\rbr{s},
\]
and
\[
	\int_{0}^{\tilde{t}} \mathbf{1}_{s\in W} g{\rbr{\tilde{x}_s}} \rbr{\dd{\tilde{x}_{s}}}_{+} = \sum_{0 < s \le t, \Delta x_s>0} \int_{x_{s-}}^{x_s} g(z) \dd{z}.
\]
Using these and (\ref{change_of_var}) we obtain \eqref{eq:BanInd1}  by simple calculations.

In a similar way one proves \eqref{eq:BanInd1-1}.
\hfill $\square$

Now we will present an alternative proof of Theorem \ref{main}. \\

{\bf Alternative proof of Theorem \ref{main}:}
First, we will assume that $g$ is continuous.

By \cite[Proposition 2.9 and Lemma 3.4]{LochowskiGhomrasni:2014},
for any $c>0$ there exists a function $x^{c}\in V^{0}$ with the
following properties: for any $s\in [0, + \ns)$,
\begin{equation}
\left|x(s)-x^{c}(s)\right|\le\frac{c}{2},\quad\left|\Delta x_{s}^{c}\right|\le\left|\Delta x_{s}\right|, \label{eq:wl1}
\end{equation}
\begin{equation}
\TTV{x^{c}}{[0,s]}{} \le\TTV x{[0,t]}{}+c \label{eq:wl2}
\end{equation}
and
\[
\Ucross{z}{x^c}{[0, s]} = \Ucross{z,c}{x}{[0, s]}.
\]
Using Lemma \ref{Ban_Viatli_extended}, for $x^{c}$ we have
\begin{align}
\int_{\R} \Ucross{z}{x^c}{[0,t]}g(z)\dd{z=} & \int_{\R}\Ucross{z,c}{x}{[0, t]}g(z)\dd z\nonumber \\
= & \int_{(0,t]}g\rbr{x_{s}^{c}}\rbr{\dd{x_{s}}^{c,cont}}_{+}+\sum_{0<s\le t,\Delta x_{s}^{c}>0}\int_{x_{s-}^{c}}^{x_{s}^{c}}g(z)\dd z,\label{eq:approx_xc}
\end{align}
where $x^{c,cont}$ denotes the continuous part of $x^c$.

By the Lebesgue dominated convergence theorem and by the definition
of $\Ucross{z}{x}{[0, t]}$, we get
\[
\lim_{c\ra0+}\int_{\R}\Ucross{z,c}{x}{[0, t]}g(z)\dd z=\int_{\R}\Ucross{z}{x}{[0, t]}g(z)\dd z.
\]
Also, to deal with the last expression in (\ref{eq:approx_xc}), for arbitrary $\varepsilon>0$
we can approximate $x$ by a step function $\tilde{x}=\sum a_{i}{\bf 1}_{[t_{i},t_{i+1})}$
such that $\sup_{s\in[0,t]}\left|g\rbr{x_{s}}-g\rbr{\tilde{x}_{s}}\right|\le\varepsilon/2$.
From (\ref{eq:wl2}) we know that
\[
\lim_{c\ra0+}\left|\int_{(0,t]}g\rbr{x_{s}^{c}}\rbr{\dd{x_{s}}^{c,cont}}_{+}-\int_{(0,t]}g\rbr{\tilde{x}_{s}}\rbr{\dd{x_{s}}^{c,cont}}_{+}\right|\le\frac{\varepsilon}{2}\TTV x{[0,t]}{}
\]
Also, by (\ref{eq:wl2}) and by the fact that $\tilde{x}$ is a step
function, we have
\[
\lim_{c\ra0+}\left|\int_{(0,t]}g\rbr{\tilde{x}_{s}}\rbr{\dd{x_{s}}^{c,cont}}_{+}-\int_{(0,t]}g\rbr{\tilde{x}_{s}}\rbr{\dd{x_{s}}^{cont}}_{+}\right|=0.
\]
Finally
\[
\lim_{c\ra0+}\left|\int_{(0,t]}g\rbr{x_{s}}\rbr{\dd{x_{s}}^{cont}}_{+}-\int_{(0,t]}g\rbr{\tilde{x}_{s}}\rbr{\dd{x_{s}}^{cont}}_{+}\right|\le\frac{\varepsilon}{2}\TTV x{[0,t]}{}.
\]
Since $\varepsilon$ may be arbitrary close to $0$, the last three relations imply
\[
\lim_{c\ra0+}\int_{(0,t]}g\rbr{x_{s}^{c}}\rbr{\dd{x_{s}}^{c,cont}}_{+}=\int_{(0,t]}g\rbr{x_{s}}\rbr{\dd{x_{s}}^{cont}}_{+}.
\]

Similarly, dividing for example the jumps of $x$ into a set of small
jumps (smaller than some threshold) whose sum of magnitudes is close
to $0$, and a finite number of bigger jumps (larger or equal than
some threshold), we get that
\[
\lim_{c\ra0+}\sum_{0<s\le t,\Delta x_{s}^{c}>0}\int_{x_{s-}^{c}}^{x_{s}^{c}}g(z)\dd z=\sum_{0<s\le t,\Delta x_{s}>0}\int_{x_{s-}}^{x_{s}}g(z)\dd z.
\]

Thus, we get that the relation (\ref{eq:BanInd1}) holds also for
any continuous $g:\R\ra\R$ and any $x\in V.$

Now we are going to drop the assumption about continuity of $g$. Let us consider the measures assigning to the Borel set $A$
the value
\[
\mu^{x}(A):=\int_{x^{-1}(A)\cap(0,t]}\rbr{\dd{x_{s}}^{cont}}_{+}.
\]
Let us notice that for any Borel measurable and bounded $g:\R\ra\R$,
\[
\int_{(0,t]}g\rbr{x_{s}}\rbr{\dd{x_{s}}^{cont}}_{+}=\int_{\R}g(y)\dd{\mu^{x}(y)}.
\]

Let ${\cal B}$ be the Borel $\sigma$-field of subsets of real numbers and let us consider the
measure $\nu^{x}$ such that
\[
\dd{\nu^{x}}(z)=\Ucross{z}{x}{[0, t]} \dd z+\dd{\mu^{x}}(z)+\sum_{0<s\le t,\Delta x_{s}>0}{\bf 1}_{[x_{s-},x_{s}]}(z)\dd z.
\]
Let $\delta>0$. By Luzin's theorem, see \cite[Theorem 7.10]{Folland1999}, for any $\varepsilon>0$ there
exists a continuous $\tilde{g}:\R\ra\R$ and $E\in{\cal B}$ such
that $\nu^{x}(E^{c})<\varepsilon$ and $\tilde{g}=g$ on $E$. Since
$g$ is bounded, taking $\varepsilon$ sufficiently small we get that
\[
\left|\int_{\R}\Ucross{z}{x}{[0, t]} g(z)\dd z-\int_{\R}\Ucross{z}{x}{[0, t]}  \tilde{g}(z)\dd z\right|<\delta
\]
\[
\left|\int_{\R}g(y)\dd{\mu^{x}(y)}-\int_{\R}\tilde{g}(y)\dd{\mu^{x}(y)}\right|<\delta
\]
and
\[
\left|\sum_{0<s\le t,\Delta x_{s}>0}\int_{x_{s-}}^{x_{s}}g(z)\dd z-\sum_{0<s\le t,\Delta x_{s}>0}\int_{x_{s-}}^{x_{s}}\tilde{g}(z)\dd z\right|<\delta.
\]
Relation (\ref{eq:BanInd1}) holds for $\tilde{g}$, since it is continuous,
thus, by  arbitrarity of $\delta$ we get that it must hold for $g$
as well.
\hfill $\square$

\section{Change of variable formulas}

In this section, we extend formula \eqref{ItoMeyyy} in two ways. The first extension assumes that $x \in V^0$ and $g: \R \ra \R$ is Borel-measurable and locally integrable (with respect to the Lebesgue measure). Under these assumptions we have:

\begin{prop} \label{prop1} If $x \in V^0$, $g: \R \ra \R$ is Borel-measurable and locally integrable (with respect to the Lebesgue measure), and one defines $f: \R \ra \R$ by $f(y) := f_0 + \int_{[0,y]} g(z) \dd z$ if $y \ge 0 $, $f(y) := f_0 - \int_{[y,0]} g(z) \dd z$ if $y < 0 $ with $f_0 \in \R$, or, equivalently, if the function $f: \R \ra \R$  is absolutely continuous and  $f(y) = f(0) + \int_{[0,y]} g(z) \dd z= f(0) + \int_{[0,y]} f'(z) \dd z$ if $y \ge 0 $ and $f(y) = f(0) - \int_{[y,0]} g(z) \dd z= f(0) - \int_{[0,y]} f'(z) \dd z$ if $y < 0 $, for some  Borel-measurable and locally integrable $g: \R \ra \R$,  then
\begin{align}
f\rbr{x_t} - f\rbr{x_0} &=   \int_{\R}\cbr{\Ucrossemph{z}{x}{\sbr{0,t}} -\Dcrossemph{z}{x}{\sbr{0, t}} }  g(z)\dd z  \nonumber \\
& =\int_{\left(0,t\right]}g\rbr{x_{s-}}{\dd{x_{s}}}+\sum_{0<s\le t,\Delta x_{s}\neq 0} \cbr{f\rbr{x_{s}}-f\rbr{x_{s-}} - g\rbr{x_{s-}}\Delta x_s}\nonumber \\
 & =\int_{\left(0,t\right]}g\rbr{x_{s}}{\dd{x_{s}}^{cont}}+\sum_{0<s\le t,\Delta x_{s}\neq 0}\cbr{f\rbr{x_{s}}-f\rbr{x_{s-}} } \nonumber \\
& = \int_{\R} \ell^z(t) g(z)\dd z + \sum_{0<s\le t,\Delta x_{s}\neq 0}\cbr{f\rbr{x_{s}}-f\rbr{x_{s-}} }. \label{TM1}
 \end{align}
 \end{prop}
{\bf Proof:} Let us apply the same notation as in the proof of Lemma \ref{Ban_Viatli_extended}. For $i \in \cbr{1,2,\ldots, N}$ the function $x$ is monotonic on $\left[t_i, t_{i+1}\right)$. This gives
\[
 \int_{\R}\cbr{\Ucross{z}{x}{\sbr{t_i, t_{i+1}}} -\Dcross{z}{x}{\sbr{t_i, t_{i+1}}} }  g(z)\dd z  = f\rbr{x_{t_{i+1}}} - f\rbr{x_{t_{i}}}.
\]
Also, for $z \notin R$, we have
\[
\sum_{i=1}^N \Ucross{z}{x}{\sbr{t_i, t_{i+1}}} = \Ucross{z}{x}{\sbr{0, t}}, \quad  \sum_{i=1}^N \Dcross{z}{x}{\sbr{t_i, t_{i+1}}}  = \Dcross{z}{x}{\sbr{0, t}}.
\]
Moreover, by \eqref{eq:BanInd1-11}
\begin{align*}
& \int_{\R}\cbr{\Ucross{z}{x}{\sbr{t_i, t_{i+1}}} -\Dcross{z}{x}{\sbr{t_i, t_{i+1}}} }  g(z)\dd z \nonumber \\
& =\int_{\left(t_i,t_{i+1}\right]}g\rbr{x_{s-}}{\dd{x_{s}}}+\sum_{t_i<s\le t_{i+1},\Delta x_{s}\neq 0}\int_{x_{s-}}^{x_{s}}\cbr{g(z)-g\rbr{x_{s-}}}\dd z\nonumber \\
 & =\int_{\left(t_i,t_{i+1}\right]}g\rbr{x_{s}}{\dd{x_{s}}^{cont}}+\sum_{t_i<s\le t_{i+1},\Delta x_{s}\neq 0}\int_{x_{s-}}^{x_{s}}g(z)\dd z.
\end{align*}
Summing the last display over $i=1,2,\ldots, N$ and applying two previous displays, we get
\begin{align*}
f\rbr{x_t} - f\rbr{x_0}
&=   \int_{\R}\cbr{\Ucross{z}{x}{\sbr{0,t}} -\Dcross{z}{x}{\sbr{0, t}} }  g(z)\dd z \nonumber \\
& =\int_{\left(0,t\right]}g\rbr{x_{s-}}{\dd{x_{s}}}+\sum_{0<s\le t,\Delta x_{s}\neq 0} \cbr{f\rbr{x_{s}}-f\rbr{x_{s-}} - g\rbr{x_{s-}}\Delta x_s}\nonumber \\
 & =\int_{\left(0,t\right]}g\rbr{x_{s}}{\dd{x_{s}}^{cont}}+\sum_{0<s\le t,\Delta x_{s}\neq 0}\cbr{f\rbr{x_{s}}-f\rbr{x_{s-}} }.
 \end{align*}
From \eqref{fiiirst} and the analogous relations for $\text{{Card}} \rbr{ {\cal D}(z) \cap (0, t]}$ and $\Dcross{y}{x}{[s,t]}$, we also infer
\begin{align*}
& \int_{\R}\cbr{\Ucross{z}{x}{\sbr{0,t}} -\Dcross{z}{x}{\sbr{0, t}} }  g(z)\dd z  \nonumber \\
& = \int_{\R} \ell^z(t) g(z)\dd z + \sum_{0<s\le t,\Delta x_{s}\neq 0}\int_{x_{s-}}^{x_{s}}g(z)\dd z \\
& = \int_{\R} \ell^z(t) g(z)\dd z + \sum_{0<s\le t,\Delta x_{s}\neq 0}\cbr{f\rbr{x_{s}}-f\rbr{x_{s-}} }.
 \end{align*}
Last two displays give \eqref{TM1}.
\hfill $\square$ \\

In the next proposition, we relax the assumption that $x \in V^0$ but then we need to assume that $g$ is locally  bounded. In the end of this section, we give example that this assumption is necessary. (Interestingly, it is possible to go into different direction -- to assume that $x$ and $f$ are c\`adl\`ag and non-decreasing -- and obtain a relevant change of variable formula, see \cite[Theorem 2.1]{BertoinYor:2013}).

\begin{prop} \label{prop2} If $x \in V$, $g: \R \ra \R$ is Borel-measurable and locally  bounded, and one defines $f: \R \ra \R$ by $f(y) := f_0 + \int_{[0,y]} g(z) \dd z$ if $y \ge 0 $, $f(y) := f_0 - \int_{[y,0]} g(z) \dd z$ if $y < 0 $ with $f_0 \in \R$, or, equivalently, if the function $f: \R \ra \R$  is locally Lipschitz and  $f(y) = f(0) + \int_{[0,y]} g(z) \dd z= f(0) + \int_{[0,y]} f'(z) \dd z$ if $y \ge 0 $ and $f(y) = f(0) - \int_{[y,0]} g(z) \dd z= f(0) - \int_{[0,y]} f'(z) \dd z$ if $y < 0 $, for some  Borel-measurable and locally integrable $g: \R \ra \R$ then
\begin{align}
f\rbr{x_t} - f\rbr{x_0} =  & \int_{\R}\cbr{\Ucrossemph{z}{x}{\sbr{0,t}} -\Dcrossemph{z}{x}{\sbr{0, t}} }  g(z)\dd z  \nonumber \\
& = \int_{\R} \ell^z(t) g(z)\dd z + \sum_{0<s\le t,\Delta x_{s}\neq 0}\cbr{f\rbr{x_{s}}-f\rbr{x_{s-}} }. \label{TM2}
 \end{align}
 \end{prop}
{\bf Proof:} Let $c>0$ and let us apply the same notation as in the second proof of Theorem \ref{main}. Using Proposition \ref{prop1} for $x^c$, \eqref{eq:approx_xc} we get
 \begin{align*}
f\rbr{x^c_t} - f\rbr{x^c_0} & =\int_{\left(0,t\right]}g\rbr{x^c_{s}}{\dd{x_{s}}^{c,cont}}+\sum_{0<s\le t,\Delta x^c_{s}\neq 0}\cbr{f\rbr{x^c_{s}}-f\rbr{x^c_{s-}} } \\
& = \int_{\left(0,t\right]}g\rbr{x^c_{s}}{\dd{x_{s}}^{c,cont}}+\sum_{0<s\le t,\Delta x^c_{s}\neq 0} \int_{x^c_{s-}}^{x^c_{s}} g(z) \dd z .
\end{align*}
Assuming that $g$ is continuous and proceeding with $c>0$ to $0$ similarly  as in the second proof of Theorem \ref{main} we get
\begin{align*}
f\rbr{x_t} - f\rbr{x_0} &=\int_{\left(0,t\right]}g\rbr{x_{s}}{\dd{x_{s}}^{cont}}+\sum_{0<s\le t,\Delta x_{s}\neq 0} \int_{x_{s-}}^{x_{s}} g(z) \dd z.
\end{align*}
Now, applying Theorem \ref{main}  we get
\begin{align*}
f\rbr{x_t} - f\rbr{x_0} &=\int_{\left(0,t\right]}g\rbr{x_{s}}{\dd{x_{s}}^{cont}}+\sum_{0<s\le t,\Delta x_{s}\neq 0} \int_{x_{s-}}^{x_{s}} g(z) \dd z \\
&= \int_{\R}\cbr{\Ucross{z}{x}{[0,t]} -\Dcross{z}{x}{[0,t]} }  g(z)\dd z.
\end{align*}
From \eqref{fiiirst} and the analogous relations for $\text{{Card}} \rbr{ {\cal D}(z) \cap (0, t]}$ and $\Dcross{y}{x}{[s,t]}$, we also infer
\begin{align*}
& \int_{\R}\cbr{\Ucross{z}{x}{\sbr{0,t}} -\Dcross{z}{x}{\sbr{0, t}} }  g(z)\dd z \\
& = \int_{\R} \ell^z(t) g(z)\dd z + \sum_{0<s\le t,\Delta x_{s}\neq 0}\int_{x_{s-}}^{x_{s}}g(z)\dd z \\
& = \int_{\R} \ell^z(t) g(z)\dd z + \sum_{0<s\le t,\Delta x_{s}\neq 0}\cbr{f\rbr{x_{s}}-f\rbr{x_{s-}} }.
 \end{align*}
Using the last two displays, we get \eqref{TM2} for continuous $g$ .

To prove \eqref{TM2} for any  bounded but not necessarily continuous $g$, we apply a mollifier function $h: \R \ra [0, +\ns)$, such that $h$ is infinitely many times differentiable, $h(y) = 0$ if $|y|>1$ and $\int_{\R} h(y) \dd y = 1$. Define
\[
g_n(z) = n \int_{\R} h(ny)g(z-y) \dd y = n \int_{\R} h(n(z-y))g(y) \dd y, \quad n\in \N.
\]
$g_n$ is continuous and bounded on any interval $[-K, K]$, $K>0$ (by the supremum of $|g|$ on $[-K-1, K+1]$). Moreover,
\[
\lim_{n \ra +\ns} \int_{\R} \left| g_n(z)- g(z) \right| \dd z = 0
\]
(see \cite[Theorem 8.14(a)]{Folland1999}). The Lebesgue dominated convergence theorem together with the local  boundedness of $\left| g_n\right|$,   $\left| g\right|$ and $\left| g_n- g \right|$ gives
\[
 \int_{\R} \cbr{\Ucross{z}{x}{\sbr{0,t}} -\Dcross{z}{x}{\sbr{0, t}} }   g_n(z)\dd z \ra \int_{\R} \cbr{\Ucross{z}{x}{\sbr{0,t}} -\Dcross{z}{x}{\sbr{0, t}} }   g(z)\dd z,
 \]
\[
 \int_{\R} \ell^z(t) g_n(z)\dd z \ra \int_{\R} \ell^z(t) g(z)\dd z,
 \]
 \[
 \sum_{0<s\le t,\Delta x_{s}\neq 0} \left|  \int_{x_{s-}}^{x_{s}}  \left| g_n(z) - g(z) \right| \dd z\right| \ra 0.
 \]
The last convergence is easily obtained when we divide the jumps of $x$ into two groups - finite number of large jumps and jumps, whose sum of magnitudes is less than a given $\varepsilon >0$. Now, for any $y \in \R$, defining $f_n(y) = f_0 + \text{sgn}(y) \int_{\sbr{0 \wedge y, 0 \vee y}} g_n(z) \dd z$
\[
f_n(y) \ra  f_0 + \text{sgn}(y) \int_{\sbr{0 \wedge y, 0 \vee y}} g(z) \dd z = f(y)
\]
as $n \ra +\ns$.
Using these limits and applying formula \eqref{TM2} for $f_n$ (and $g_n$) we obtain \eqref{TM2} for $f$  and $g$.
\hfill $\square$ \\

We end this section with an easy example that the assumption that $f$ is locally Lipschitz in Proposition \ref{prop2} is necessary.
Let
\[
x(t) = \begin{cases}
\frac{1+(-1)^{\lfloor \frac{1}{1-t}\rfloor^2}}{2\lfloor \frac{1}{1-t}\rfloor^2} \text{ if } t \in [0,1), \\
0 \text{ if } t \ge 1;
\end{cases}
\quad
f(y) = \begin{cases}
\sqrt{y} \text{ if } y \in [0,+\ns), \\
0 \text{ if } y < 0.
\end{cases}
\]
$f$ is absolutely continuous. $x$ has finite variation, it has no continuous part and has two jumps, from $1/(2n)^2$ to $0$ and from $0$ to $1/(2n)^2$, for each $n \in \N$, on the interval $t \in (0,1)$. From this it follows that the sum $$\sum_{0 < s \le 1, \Delta x_s \neq 0} \cbr{f\rbr{x_s} - f\rbr{x_{s-}}}$$ is not well defined.

\bibliographystyle{amsalpha}
\bibliography{/Users/rafallochowski/biblio/biblio}

\end{document}